\documentclass{svproc}
\usepackage{amsmath, verbatim, amsfonts, amssymb, color}
\usepackage{mathrsfs}
\usepackage{dsfont}

\renewcommand\labelenumi{\textup{\alph{enumi})}}
\renewcommand\theenumi\labelenumi

 \makeatletter
\def\@makefnmark{\hbox{(\@textsuperscript{\normalfont\@thefnmark})}}
\makeatother

\makeatletter
\@namedef{subjclassname@2020}{%
  \textup{2020} Mathematics Subject Classification}
\makeatother

\providecommand{\ack}[1]{\par\addvspace\baselineskip
\noindent\ackname\enspace\ignorespaces#1}%
\def\subjclassname{\textup{2020} \textit{Mathematics Subject Classification:}}%
\providecommand{\subjclass}[1]{\par\addvspace\baselineskip
\noindent\subjclassname\enspace\ignorespaces#1}%

\newcommand{\scalar}[2]{\left\langle #1 \,\middle |\, #2\right\rangle}

\newcommand{\fstop}{\,\,\mathrm{.}}

\def\M{{M}}
\def\g{{g}}
\DeclareMathOperator{\vol}{\mathsf{vol}}

\DeclareMathOperator{\Ric}{\mathrm{Ric}}
\DeclareMathOperator{\scal}{\mathsf{scal}}

\def\Leb{\mathcal{L}}
\def\C{\mathcal{C}}

\def\Hil{{H}} 
\def\N{{\mathbb N}}    
\def\Z{{\mathbb Z}} 
\def\R{{\mathbb R}}
\def\TT{\mathbb T}
\newcommand{\ee}{e}
\begin{document}
\mainmatter
\title{\bfseries Random Riemannian Geometry in 4 Dimensions}
\titlerunning{Random Riemannian Geometry in 4 Dimensions}

\author{Karl-Theodor Sturm}
\authorrunning{Sturm}
\institute{Hausdorff Center for Mathematics, University of Bonn, Germany \\
\email{sturm@uni-bonn.de}  }

\maketitle

\begin{abstract}  
We construct and analyze conformally invariant random fields on 4-dimensional Riemannian manifolds $(M,g)$. These centered Gaussian fields $h$, called \emph{co-biharmonic Gaussian fields}, are characterized by their covariance kernels $k$ defined as the integral kernel for the inverse of the \emph{Paneitz operator}
\begin{equation*}\mathsf p=\frac1{8\pi^2}\bigg[\Delta^2+
 \mathsf{div}\left(2\mathsf{Ric}-\frac23\mathsf{scal}\right)\nabla
\bigg].
\end{equation*}
The kernel $k$  is invariant (modulo additive corrections) under conformal transformations, and it exhibits a precise logarithmic divergence $$\Big|k(x,y)-\log\frac1{d(x,y)}\Big|\le C.$$
In terms of the co-biharmonic Gaussian  field $h$, we define the \emph{quantum Liouville measure}, a random measure on $M$,
heuristically given as
\begin{equation*}
  d\mu(x):= e^{\gamma h(x)-\frac{\gamma^2}2k(x,x)}\,d \vol_\g(x)\,,
\end{equation*}
and rigorously obtained a.s.~for $|\gamma|<\sqrt8$ as weak limit of the RHS with $h$ replaced by suitable regular approximations $(h_\ell)_{\ell\in\N}$.

For the flat torus $M=\TT^4$, we provide discrete approximations of the Gaussian field and of the Liouville measures in terms of semi-discrete random objects, based on Gaussian random variables on the discrete torus and piecewise constant functions in the isotropic  Haar system.
\ack{ The author gratefully acknowledges financial support from the European Research Council through the ERC AdG `RicciBounds' (grant agreement 694405)
as well as funding by the Deutsche Forschungsgemeinschaft through the project `Random Riemannian Geometry' within the SPP 2265 `Random Geometric Systems'. The author also likes to thank the reviewer for his/her  careful reading and valuable comments.}
\end{abstract}

\bigskip\noindent

\section{Random Riemannian Geometries and Conformal Invariance}

The basic ingredients of any \emph{Random Riemannian Geometry} are a family $\mathfrak M$ of Riemannian manifolds $(M,g)$ and a probability measure $\mathbf P_{\mathfrak M}$ on $\mathfrak M$. 
Typically, $\mathfrak M=\{(M,g'): \ g'= \ee^{2h}g, \ h\in\C^\infty(M)\}$ for some given $(M,g)$, and $\mathbf P_{\mathfrak M}$ is the push forward of a probability measure $\mathbf P_{g}$ on $\mathcal C^\infty(M)$ under the map $h\mapsto (M, e^{2h}g)$.

Of major interest are Random Riemannian Geometries which are \emph{conformally invariant}. In the previous setting this means that 
\begin{itemize}
\item \ $\mathbf P_{g'}=\mathbf P_{g}$ if $g'= \ee^{2\varphi}g$ for some $\varphi\in\C^\infty(M)$
\item \ $h\stackrel{(d)}=h'\circ \Phi$ if  $\Phi: (M,g)\to (M',g')$ is an isometry and $h$ and $h'$ are distributed according to $\mathbf P_{g}$ and  $\mathbf P_{g'}$, resp.
\end{itemize}
In this respect, of course,   the 2-dimensional case 
plays a particular role   thanks to the powerful Riemannian Mapping Theorem 
--- but the concept of conformally invariant random geometries is by no means restricted to this case.

Mostly, such probability measures $\mathbf P_{g}$ are Gaussian fields, informally given as
\begin{equation}\label{gauss}d\mathbf P_{g}(h)=\frac1{Z_{g}}\,\exp\Big(- \mathfrak e_g (h,h)\Big)\,dh\end{equation}
with some (non-existing) uniform distribution $dh$ on  $\C^\infty(M)$, a normalizing constant $Z_g$, and some bilinear form $\mathfrak e_g$. 
The rigorous definition of such probability measures $\mathbf P_{g}$ often requires to pass to spaces of distributions (rather than smooth functions). It is based on the Bochner--Minlos Theorem and the unique characterization of $\mathbf P_{g}$ as
\begin{equation} \int e^{i\langle u,h\rangle}\,d \mathbf P_{g}(h)=\exp\Big(-\frac12 \mathfrak k_g(u,u)\Big)\qquad\quad\forall u\in\C^\infty(M)
\end{equation}
where $\mathfrak k_g(u,u)^{1/2}:=\sup_{h}\frac{\langle u,h\rangle}{\mathfrak e_g(h,h)^{1/2}}$
 denotes the norm 
dual to  $\mathfrak e_g$.

%
%

\medskip
The conformal invariance requirement for the random geometry then amounts to the requirement
\begin{equation}\label{conf-inv-en}
\mathfrak e_g(u,u)=\mathfrak e_{e^{2\varphi}g}(u,u)\qquad\quad\forall \varphi,\, \forall u.\end{equation}
In the two-dimensional case, this is a well-known property of the \emph{Dirichlet energy}, cf.~\cite{FOT},
$$\mathcal E_g(u,u):=\int_M \big|\nabla_g u\big|^2\,d\vol_g.$$ 

\medskip
The conformally invariant random field defined and constructed in this way is  the celebrated \emph{Gaussian Free Field}~\cite{She07}. It 
is a particular (and  the most prominent) case of a log-correlated random field~\cite{DRSV} and of a fractional Gaussian field~\cite{LSSW}.
It naturally arises as the scaling limit of various discrete models of random surfaces, for instance discrete Gaussian Free Fields  or harmonic crystals~\cite{She07}.
It is also deeply related to another planar conformally invariant random object of fundamental importance, the \emph{Schramm--Loewner evolution} \cite{Lawler05,Lawler18,Schramm07}.
For instance, level curves of the Discrete Gaussian Free Field  converge to~$\mathsf{SLE}_{4}$~\cite{SchShe09}, and zero contour lines of the Gaussian Free Field are well-defined random curves distributed according to~$\mathsf{SLE}_{4}$~\cite{SchShe13}.
The work~\cite{MSImaginary1}, and subsequent works in its series, thoroughly study the relation between the Schramm-Loewner evolution and Gaussian free field on the plane.
The \emph{Liouville Quantum Gravity} is a random measure, informally obtained as the Riemannian volume measure when the metric tensor is conformally transformed with the Gaussian Free Field as conformal weight.
Since the Gaussian Free Field is only a distribution, the rigorous construction of  the random measure requires a  renormalization procedure due to Kahane \cite{Kah85}.
This renormalization depends on a roughness parameter $\gamma$ and works only for $|\gamma| < 2$.
In~\cite{MSLQG1} and subsequent work in its series, Miller and Sheffield prove that for the value $\gamma = \sqrt{8/3}$ the Liouville Quantum Gravity coincides with the Brownian map, that is a random metric measure space arising as a universal scaling limit of random trees and random planar graphs (see \cite{LeGall19,LeGallMiermont12} and the references therein).
More recently, \cite{DDDF20,GM21} establish the existence of the Liouville Quantum Gravity metric for $\gamma \in (0,2)$.

\medskip
All these approaches to conformally invariant random objects so far, with exception of the recent contribution \cite{Cercle}, are limited to the two-dimensional case. The main reason is not the lack of a Riemannian Mapping Theorem but the fact that the  Dirichlet energy is  no longer conformally invariant  in dimension $n\not=2$. One rather obtains
$$\mathcal E_{e^{2\varphi}g}(u,u)=\int_M \big|\nabla_g u\big|^2\,\ee^{(n-2)\varphi}d\vol_g.$$
In the four-dimensional case, a more promising candidate appears to be the bi-Laplacian energy
$$\tilde{\mathfrak e}_g(u,u):=\int_M \big(\Delta_g u\big)^2\,d\vol_g.$$ 
This energy functional is still not conformally invariant but it is close to: 
$$\tilde{\mathfrak e}_{e^{2\varphi}g}(u,u)=\int_M \big(\Delta_g u+2\nabla_g\varphi\, \nabla_gu\big)^2\,d\vol_g=\tilde{\mathfrak e}_g(u,u)+\text{ low order terms}.$$ 
Our search for a conformally invariant energy functional in dimensions $n=4$ finally will lead us to considering
$$\mathfrak e_g(u,u)=c\int_M \big(\Delta_g u\big)^2\,d\vol_g+\text{ low order terms}.$$ 
Paneitz \cite{Paneitz} found the  precise formula for the conformally invariant energy functional in dimension 4. Subsequently, Graham, Jenne, Mason, and Sparling \cite{GJMS} showed the existence of a conformally invariant energy functional of the form
$$\mathfrak e_g(u,u)=c\int_M \big(-\Delta_g u\big)^{n/2}\,d\vol_g+\text{ low order terms}$$ 
on Riemannian manifolds of even dimension $n$.

\medskip
Based on these results, jointly with L.~Dello Schiavo, R.~Herry and E.~Kopfer \cite{DHKS}, we   constructed and analyzed conformally invariant random fields on Riemannian manifolds $(M,g)$ of arbitrary even dimensions. In the subsequent Sections 2 and 3, we will summarize these results in the particular case $n=4$. In Section 4, we will present a detailed study of approximations of the random field and the random measure on the 4-dimensional flat torus in terms  of corresponding random objects on the discrete tori $\TT^4_\ell$, $\ell\in\N$.

\section{Paneitz Energy on 4-Dimensional Manifolds} 
From now, we will be more specific. $(M,g)$ will always  be a 4-dimensional smooth, compact, connected Riemannian manifold without boundary.
Integrable functions (or distributions) $u$ on $M$ will be called \emph{grounded} if $\int_M u\,d\vol_g=0$ (or $\langle u,\bf 1\rangle=0$, resp.).
Let $(\varphi_j)_{j\in\N_0}$ denote the complete ON-basis of $L^2(M,\vol_g)$ consisting of eigenfunctions of $-\Delta_g$ with corresponding eigenvalues $(\lambda_j)_{j\in\N_0}$. 
Then the \emph{grounded Sobolev spaces} $\mathring H^s(M,g)=(-\Delta_g)^{-s/2} \mathring L^2(M,\vol_g)$ for $s\in\R$ are given
by
$$\mathring H^s(M,g)=\bigg\{ u=\sum_{j\in\N} \alpha_j\varphi_j: \ \sum_{j\in\N}\lambda_j^s\,|\alpha_j|^2<\infty \bigg\},$$
whereas the usual Sobolev spaces are $H^s(M,g)=(1-\Delta_g)^{-s/2}L^2(M,\vol_g)=\mathring H^s(M,g)\oplus\, \R\cdot\bf1$.
Extending the scalar product in $L^2(M,\vol_g)$, the pairing between $u=\sum_{j\in\N_0}\alpha_j\varphi_j\in \Hil^s$ and $v=\sum_{j\in\N_0}\beta_j\varphi_j\in \Hil^{-s}$ is given by
$$\langle u,v\rangle:=\langle u,v\rangle_{\Hil^s,\Hil^{-s}}:=\sum_{j\in\N_0} \alpha_j\,\beta_j.$$

The  Laplacian acts on these grounded spaces by $-\Delta_g: \mathring H^s\to \mathring  H^{s-2}, \ \sum_{j\in\N} \alpha_j\varphi_j\mapsto \sum_{j\in\N} \lambda_j\alpha_j\varphi_j$.
The operator inverse to it is the \emph{grounded Green operator}
$$\mathring{\sf G}_g: \mathring  H^s\to \mathring  H^{s+2}, \ \sum_{j\in\N} \alpha_j\varphi_j\mapsto \sum_{j\in\N} \frac{\alpha_j}{\lambda_j}\varphi_j.$$
On $\mathring H^0=\mathring L^2$, it is given as an integral operator
$\mathring{\sf G}_gu(x)=\int_{M} \mathring G_g(x,y) \,u(y)\,d\vol_g(y)$
in terms of the grounded Green kernel $\mathring G_g(x,y)$ on $M$.
The latter is symmetric in $x$ and $y$, it is grounded (i.e. $\int_{M}\mathring G_g(x,y) \,d\vol_g(y)=0$ for all $x$) and
$$\Big|\mathring G_g(x,y) -\frac{1}{4\pi^2\cdot  d_g(x,y)^2}\Big|\le C.$$

\begin{definition} The \emph{Paneitz energy} 
 is defined as the bilinear form on $L^2(M,\vol_g)$ with domain $H^2(M)$ by
\begin{equation}
\mathfrak e_g(u,u)=\frac1{8\pi^2}\int_M \bigg[ (\Delta_g u)^2 - 2\Ric_g(\nabla_g u,\nabla_gu)+\frac23\scal_g\cdot|\nabla_gu|^2 \bigg]\,d\vol_g.
\end{equation}
\end{definition}
In particular, for every 4-dimensional Einstein manifold with $\Ric_g=k\,g$ for $k\in\R$ (which implies $\scal_g=4k$),
\begin{equation}
\mathfrak e_g(u,u)=\frac1{8\pi^2}\int_M \bigg[ (\Delta_g u)^2 +\frac23k\, |\nabla_gu|^2 \bigg]\,d\vol_g.
\end{equation}

\begin{example} a) For the 4-sphere $M={\mathbb S}^4$,
\begin{equation}
\mathfrak e_g(u,u)=\frac1{8\pi^2}\int_M \bigg[ (\Delta_g u)^2 +2 |\nabla_gu|^2 \bigg]\,d\vol_g.
\end{equation}
b) For the 4-torus $M=\TT^4$,
\begin{equation}
\mathfrak e_g(u,u)=\frac1{8\pi^2}\int_M (\Delta_g u)^2 \,d\vol_g.
\end{equation}
\end{example}

\begin{theorem}[\cite{Paneitz}] The Paneitz energy is conformally invariant:
$$
\mathfrak e_g(u,u)=\mathfrak e_{e^{2\varphi}g}(u,u)\qquad\quad\forall \varphi\in\C^\infty(M),\, \forall u\in H^2(M).$$
\end{theorem}

\begin{definition} The 4-manifold $(M,g)$ is called \emph{admissible} if $\mathfrak e_g>0$ on $\mathring H^2(M)$.
\end{definition}
As an immediate consequence of Theorem 1, we observe that admissibility is a conformal invariance. Large classes of 4-manifolds are admissible.
\begin{proposition}[{\cite[Prop.~2.4, 2.5]{DHKS}}]
a) All compact Einstein 4-manifolds with nonnegative Ricci curvature 
are admissible.

b) All compact hyperbolic 4-manifolds with spectral gap $\lambda_1>2$ are admissible.
\end{proposition}
However, not every compact four-dimensional Riemannian manifold is admissible.

\begin{example}[{\cite[Prop.~2.7]{DHKS}}] \label{non-admiss}
Let $M_1, M_2$ be compact hyperbolic Riemannian surfaces such that $\lambda_1(M_1)\le\frac23$. 
Then the Einstein 4-manifold $M=M_1\times M_2$ is not admissible. 
 \end{example}

If $(M,g)$ is admissible, then the \emph{Paneitz operator} (or \emph{co-bilaplacian})
\begin{equation}\label{pan-op}
\mathsf p_g=\frac1{8\pi^2}\bigg[\Delta_g^2+
 \mathsf{div}\left(2\mathsf{Ric}_g-\frac23\mathsf{scal}_g\right)\nabla
\bigg]
\end{equation}
is a self-adjoint positive operator on $L^2(\M,\vol_g)$ with  domain $\Hil^4(M)$.
Here  the curvature term $2\mathsf{Ric}_g-\frac{2}{3}\mathsf{scal}_g$ should be viewed as an endomorphism of the
tangent bundle, acting on the gradient of a function.
In coordinates: 
\begin{equation*}
  \mathsf{p}_g u=\frac1{8\pi^2}\,\sum_{i,j}\nabla_i\left[\nabla^i\nabla^j+2 \mathsf{Ric}^{ij}_g-\frac23\mathsf{scal}_g\cdot g^{ij}\right] \nabla_j u, \qquad \forall u \in \mathcal{C}^{\infty}(M).
\end{equation*}

Let $(\psi_j)_{j\in\N_0}$ denote a complete orthonormal basis of $L^2(M,\vol_g)$ consisting of eigenfunctions for $\mathsf p_g$, and let $(\nu_j)_{j\in\N_0}$ denote the corresponding sequence of eigenvalues.
Then the operator $\mathsf k_g$, inverse to $\mathsf p_g$ on $\mathring L^2$,  is given on $\Hil^{-4}(M)$ by
$$\mathsf k_g: \ u\mapsto {\mathsf k}_g u:=\sum_{j\in\N}\frac1{\nu_j}\, \langle u,\psi_j\rangle\, \psi_j,$$
and the associated bilinear form with domain $\Hil^{-2}(M)$ is given by
$${\mathfrak k}_g(u,v):=\langle u,\mathsf k_g v\rangle_{L^2}=\sum_{j\in\N}\frac1{\nu_j}\, \langle u,\psi_j\rangle\, \langle v,\psi_j\rangle.$$
The  crucial properties of the kernel for the \emph{co-biharmonic Green operator} ${\sf k}_g$ are its  logarithmic divergence and its conformal invariance.

\begin{theorem}[{\cite[Thm.~2.18]{DHKS}}]  If $(M,g)$ is admissible, then $\mathsf k_g$ is an integral operator with an integral kernel $k_g$  which satisfies
\begin{equation}\label{log-kernel}
\Big|k_g(x,y) +\log d_g(x,y)\Big|\le C.
\end{equation}
Furthermore, the kernel $k_g(x,y)$ is symmetric in $x,y$ and grounded.
\end{theorem}
\begin{theorem}[{\cite[Prop.~2.19]{DHKS}}] \label{conf-k}\label{t:Covariant} 
  Assume that $(M,g)$ is admissible and that $g':= e^{2\varphi}g$ for some $\varphi\in\mathcal{C}^\infty(M)$.
  Then the co-biharmonic Green kernel ${k}_{g'}$ for the metric $g'$ 
is given by 
\begin{align}\label{trafo-k}
k_{g'}(x,y)=& \ k_g(x,y)-\frac12\bar\phi(x) -\frac12\bar\phi(y)
\end{align}
with $\bar\phi\in\C^\infty(M)$  defined by 
\begin{align*}\bar\phi&:= \frac2{\mathsf{vol_{g'}}(M)}\int k_g(.,z)\,d\mathsf{vol}_{g'}(z)-\frac1{\mathsf{vol_{g'}}(M)^2}\iint k_g(z,w)\,d\mathsf{vol}_{g'}(z)\,d\mathsf{vol}_{g'}(w)\fstop
 \end{align*}\end{theorem}
\begin{example} Assume that $(M,g)$ is Ricci flat. Then 
$$k_g(x,y)=8\pi^2\,\mathring G_g^{(2)}(x,y):=8\pi^2\,\int_M \mathring G_g(x,z)\,\mathring G_g(z,y)\,d\vol_g(z)$$
where $\mathring G_g$ denotes the grounded Green kernel on $(M,g)$.
\end{example}
\section{Co-biharmonic Gaussian Field and Quantum Liouville Measure}
Throughout the sequel, assume that $(M,g)$ is an admissible 4-manifold (compact, smooth, without boundary -- as always).
\subsection{Conformally Invariant Gaussian Field}
\begin{definition}
A co-biharmonic Gaussian field $h$ on~$(M,g)$ is a linear family $$\big(\langle h, u\rangle\big)_{u\in {H}^{-2}}$$ of centered Gaussian random variables (defined on some probability space) 
with
$$\mathbf{E}\big[\langle h, u\rangle^2\big]={\mathfrak k}_g(u,u)\qquad\quad\forall u\in H^{-2}(M).$$
\end{definition}

\begin{theorem}[{\cite[Prop.~3.9, Rem.~3.3]{DHKS}}]  Let a probability space $(\Omega, \mathfrak F, \mathbf P)$ be given and an i.i.d. sequence $(\xi_{j})_{j\in\N}$ of ${\mathcal N}(0,1)$ random variables. Furthermore, let $(\psi_j)_{j\in\N_0}$ and $(\nu_j)_{j\in\N_0}$ denote the sequences of eigenfunctions and eigenvalues for $\mathsf p_g$ (counted with multiplicities). Then a co-biharmonic field is given by
\begin{equation}
h:=\sum_{j\in\N} \nu_j^{-1/2}\,\xi_j\,\psi_j. \end{equation}
More precisely,
\begin{enumerate}
\item For each $\ell\in\N$, a centered Gaussian random variable $h_\ell$ with values in $\C^\infty(M)$ is given by
\begin{equation}
h_\ell:=\sum_{j=1}^\ell \nu_j^{-1/2}\,\xi_j\,\psi_j.\end{equation}
\item The convergence $h_\ell\to h$ holds in $L^2(\mathbf P)\times \Hil^{-\epsilon}(M)$ for every $\epsilon>0$. In particular,
for a.e.~$\omega$ and every $\epsilon>0$, 
$$h^\omega\in \Hil^{-\epsilon}(M),$$

\item For every $u\in\Hil^{-2}(M)$, the family $(\langle u,h_\ell\rangle)_{\ell\in\N}$ is a centered $L^2(\mathbf P)$-bounded martingale and
$$\langle u,h_\ell\rangle\to\langle u,h\rangle\quad\text{in $L^2(\mathbf P)$ as }\ell\to\infty.$$
\end{enumerate}
\end{theorem}

\begin{remark}
a) A co-biharmonic Gaussian field on $(M,g)$ can be regarded as a random variable with values in $\mathring\Hil^{-\epsilon}(M)$ for any $\epsilon>0$. 

b) Given any `grounded' white noise $\Xi$ on $(M,g)$, then $h:=\sqrt{\mathsf k_g} \Xi$ is a co-biharmonic Gaussian field on $(M,g)$.
\end{remark}
\begin{theorem}[{\cite[Thm.~3.11]{DHKS}}]  Let $h:\Omega\to H^{-\epsilon}(M)$ 
denote a co-biharmonic Gaussian field for $(M,g)$ and 
let  $\g'=e^{2\varphi}\g$ with $\varphi\in\C^\infty(\M)$.
Then
\begin{equation*}
h':=h-\frac1{\vol_{g'}(M)}
\big\langle{h,\mathbf 1\big\rangle}_{H^{-\epsilon}(M,g'),H^\epsilon(M,g')}
\end{equation*}
is a co-biharmonic Gaussian field for $(M,g')$.
\end{theorem}

Besides the previous eigenfunction approximation, there are numerous other ways to approximate a given  co-biharmonic Gaussian field $h$ by `smooth' Gaussian fields $h_\ell, \ell\in\N$. 

\begin{proposition}\label{convolution-approx}  Let $\rho_\ell$ for $\ell\in\N$  be a family of bounded functions on $\M\times\M$ such that 
$\rho_\ell(x,.)\vol_g$ for each $x\in\M$ is a family of probability measures on $\M$ which for $\ell\to\infty$ weakly converges to $\delta_x$.
Define centered Gaussian fields $h_\ell$ for $\ell\in\N$ by
\begin{equation}\label{kern-appr}
h_\ell(y):=\scalar{h}{\rho_\ell(.,y)}.
\end{equation}
Then, for every~$u\in \C(M)$,   as $\ell\to\infty$
\begin{equation*}
\scalar{h_\ell}{u}\longmapsto \scalar{h}{u}\qquad\text{
$\mathbf P$-a.s.\ and in $L^2(\mathbf P)$}\fstop
\end{equation*}
The associated covariance kernels are given by
$$k_\ell(x,y):=
\iint k(x',y')\rho_\ell(x',x)\rho_\ell(y',y)\,d\vol_\g(x')\,d\vol_\g(y')$$
for all $\ell\in\N$, and $k_\ell\to k$ as $\ell\to\infty$
 on locally uniformly on $M\times M$ off the diagonal.
\end{proposition}

\begin{proof} Obviously, 
$\scalar{h_\ell}{u}= \scalar{h}{\rho_\ell * u}$
with $(\rho_\ell*u)(x)=\int \rho_\ell(x,y) u(y) \,d\vol_\g(y)$, and $\rho_\ell* u\to u$ in $L^2$ as $\ell\to\infty$. Moreover,
\begin{align*}
{\mathbf E}\Big[\big|\langle h |  u\rangle-\langle h_\ell u\rangle\big|^2\Big]
&={\mathbf E}\Big[\big|\langle h|u-\rho_\ell*u\rangle\big|^2\Big]\\
&=\sum_{j=1}^\infty \frac1{\nu_j} \big|\langle \psi_j|u-\rho_\ell*u\rangle\big|^2
\le C \big\|u-\rho_\ell*u\big\|^2_{H^{-2}}.
\end{align*}
\end{proof}

A particular case of such approximations through convolution kernels will be considered now. 
\begin{proposition}\label{const-approx-mart}
Let $(\mathfrak Q_\ell)_{\ell\in\N}$ be a family of partitions of $\M$ with $\forall \ell, \forall Q\in \mathfrak Q_\ell: \exists m\in\N, \exists Q_1,\ldots, Q_m\in \mathfrak Q_{\ell+1}: Q=\bigcup_{i=1}^m Q_i$ and with $\sup\{\mathrm{diam}(Q): Q\in \mathfrak Q_\ell\}\to0$ as $\ell\to\infty$. For $\ell\in\N$ put
\begin{equation}\label{cons-kern}\rho_\ell(x,y):=\sum_{Q\in \mathfrak Q_\ell} \frac1{\vol_g(Q)} {\bf 1}_Q(x) {\bf 1}_Q(y).
\end{equation}
In other words, for given $x\in\M$ we have  $\rho_\ell(x,.)=\frac1{\vol_g(Q)}  {\bf 1}_Q$ with the unique ${Q\in \mathfrak Q_\ell}$ which contains $x$. Defining $h_\ell$ as before then yields

\begin{equation}\label{cons-appr}
h_\ell(x)=\frac1{\vol_g(Q)} \scalar{h}{{\bf 1}_Q}\qquad\forall x\in Q, \forall Q\in \mathfrak Q_\ell.
\end{equation}
For $\ell\in\N$, let $\mathfrak F_\ell$ denote the $\sigma$-field in $(\Omega, {\frak F},\mathbf Q)$ generated by the  random functions on $\M$ that are piecewise constant on each of the sets ${Q\in \mathfrak Q_\ell}$. Then $(h_\ell)_{\ell\in\N}$ is a $(\Omega, {\frak F},({\frak F}_\ell)_{\ell\in\N},\mathbf P)$-martingale 
and
\begin{equation}
h_\ell={\mathbf E}\big[h \big| \mathfrak F_\ell\big]\qquad\forall \ell\in\N.
\end{equation}
\end{proposition}

\subsection{Quantum Liouville Measure}

Let an admissible 4-manifold $(M,g)$ be given as well as a co-biharmonic Gaussian field $h$ on it. Furthermore, let smooth approximations $(h_\ell)_{\ell\in\N}$ of it be given --- informally defined as 
$h_\ell:=\rho_\ell*h$ and formally by \eqref{kern-appr} --- in terms of a sequence $(\rho_\ell)_{\ell\in\N_0}$ of bounded convolution densities on $M$. Fix $\gamma\in\R$.

For $\ell\in\N$ define a random measure $\mu_\ell=\rho_\ell\,\vol_g$ on $M$ with density
$$\rho_\ell(x):=\exp\Big(\gamma h_\ell(x)-\frac{\gamma^2}2k_\ell(x,x)\Big)$$
with
$k_\ell(x,y):=
\iint k(x',y') \rho_\ell(x',x)\rho_\ell(y',y)\,d\vol_\g(x')\,d\vol_\g(y')$
 as before. 
\begin{theorem}[{\cite[Thm.~4.1]{DHKS}}]  If $|\gamma|<\sqrt 8$, then there exists a random measure $\mu$ on $M$ with
$\mu_\ell\to\mu$. More precisely, for every $u\in\C(M)$,
$$\int_Mu\,d \mu_\ell \longrightarrow\int_Mu\,d \mu\quad\text{in $L^1(\mathbf P)$ and $\mathbf P$-a.s. as }\ell\to\infty.$$
The random measure $\mu$ is independent of the choice of the convolution densities $(\rho_\ell)_{\in\N}$.

If the $(\rho_\ell)_{\in\N}$ are chosen according to \eqref{cons-kern} then for each $u\in \mathcal C(M)$ the family
$Y_\ell:=\int_\M u\,d\vol_\g, \ell\in\N$, is a uniformly integrable martingale.
If in addition
$|\gamma|<2$, then this martingale is even $L^2$-bounded. 
\end{theorem}

The latter claim, indeed, can be seen directly:
 \begin{align*}
\sup_\ell{\mathbf E}\Big[{Y_\ell}^2\Big]&=\sup_\ell{\mathbf E}\iint \ee^{\gamma(h_\ell(x)+h_\ell(y)-\frac{\gamma^2}2(\mathbf E[h_\ell^2(x)+h_\ell^2(y)]} \,
u(x)u(y)\, d\text{vol}_\g(x)d\text{vol}_\g(y)\\
&=\sup_\ell\iint \ee^{\gamma^2 k_\ell(x,y)}\, u(x)u(y)\, d\text{vol}_\g(x)\,d\text{vol}_\g(y)\\
&\le \|u\|_\infty^2\cdot \sup_\ell\iint \bigg[\iint \rho_\ell(x',x)\rho_\ell(y',y)\ee^{\gamma^2  k(x',y')}\, d\text{vol}_\g(x')\,d\text{vol}_\g(y')\bigg]\\ 
&\qquad\qquad\qquad\qquad\qquad d\text{vol}_\g(x)\,d\text{vol}_\g(y)\\
&= \|u\|_\infty^2\cdot \iint\ee^{\gamma^2  k(x',y')}\, d\text{vol}_\g(x')\,d\text{vol}_\g(y')\\
&\le \|u\|_\infty^2\cdot\iint \frac1{d(x,y)^{\gamma^2}}\,d\text{vol}_\g(x)\, d\text{vol}_\g(y)+C'
\end{align*}
by means of Jensen's inequality and the kernel estimate \eqref{log-kernel}. Obviously, the final integral is finite if and only if $\gamma^2<4$.

\begin{definition} The random measure $\mu:=\lim\limits_{\ell\to\infty}\mu_\ell$ is  called \emph{quantum Liouville measure}.
\end{definition}

\begin{remark}[{\cite[Cor~4.10, Prop.~4.14]{DHKS}}] Assume $|\gamma|<\sqrt 8$ and let $\omega\mapsto\mu^\omega$ denote the random measure constructed above. Then for $\mathbf P$-a.e.~$\omega$, the measure $\mu^\omega$ on $M$
\begin{itemize}
\item 
 does not charge sets of vanishing $H^2$-capacity;
\item 
does not charge sets of vanishing $H^1$-capacity
provided $|\gamma|<2$;
\item 
is singular w.r.t.~the volume measure on $M$ whenever $\gamma\not=0$.
\end{itemize}
Moreover, the random measure $\mu$  has finite moments of any negative order, i.e.~for any $p>0$,
$${\mathbf E}\big[\mu(M)^{-p}\big]<\infty.$$
\end{remark}

A key property of the quantum Liouville measure is its quasi-invariance under conformal transformations.
\begin{theorem}[{\cite[Thm.~4.4]{DHKS}}]  Let $\mu$ be the quantum Liouville measure for $(M,g)$, and $\mu'$ be the quantum Liouville measure for $(M,g')$ where
$g'=e^{2\varphi}\g$ for some $\varphi\in\C^\infty(\M)$. Then
\begin{equation}
\mu'
\ \stackrel{{\rm (d)}}=\
e^{-\gamma\xi+ \frac{\gamma^2}{2} \bar\varphi+4\varphi}\, \mu
\end{equation}
where 
$\xi:=  \frac1{v'}\langle h,\ee^{4\varphi}\rangle$ and $\bar\varphi:= \frac2{v'}\,{\sf k}_g(e^{4\varphi})-\frac1{{v'}^2}\,{\frak k}_g(e^{4\varphi},e^{4\varphi})$ with $v':= \vol_{g'}(M)$.
\end{theorem}

\section{Approximation by Random Fields and Liouville Measures on the Discrete 4-Torus}

For the remaining part, we now focus on the 4-dimensional torus $\TT^4:=\R^4/\Z^4$, equipped with the flat metric. With this choice of $(M,g)$, we will drop the $g$ from the notations: $k=k_g, G=G_g$ etc.

For the 4-torus, we will study approximations of the co-biharmonic field --- now briefly called biharmonic field (since the underlying  Paneitz operator or co-bilaplacian  is now simply  the bilaplacian) --- and of the quantum Liouville measure by (semi-) discrete versions of such fields and measures, defined on the discrete tori $\TT^4_\ell$ as $\ell\to\infty$. 

\subsection{The Isotropic Haar System}
To begin with, for $\ell\in\N_0$ define the parameter sets
$$A_\ell:=\big\{0,1,\ldots,2^{\ell}-1\big\}^4, \quad B_\ell:=\{0,1\}^4\setminus \{(0,0,0,0)\}, \quad I_\ell:=A_\ell\times  B_\ell$$
and the
discrete 4-torus 
$$\TT^4_\ell:=2^{-\ell}\cdot A_\ell\, =\, (2^{-\ell}\Z^4)/\Z^4.$$
Moreover, let 
$\mathfrak Q_\ell:=\big\{ Q_{\ell,\alpha}: \ \alpha\in A_\ell\big\}$
denote the set of all dyadic cubes 
$$Q_{\ell,\alpha}:=2^{-\ell}\cdot\Big( [\alpha_1,\alpha_1+1)\,\times\,
 [\alpha_2,\alpha_2+1)\,\times\, [\alpha_3,\alpha_3+1)\,\times\, [\alpha_4,\alpha_4+1)\Big)
\,\subset\,\TT^4
 $$ 
of edge length $2^{-\ell}$, and let ${\mathcal S}_\ell$ denote the set of all grounded functions $u:\TT^4\to\R$ which are constant on each of the cubes $Q\in \mathfrak Q_\ell$.
With each $Q_{\ell,\alpha}\in \TT^4_\ell$ we  associate a set $\{\eta_{\ell,\alpha,\beta}: \beta\in B_\ell\}\subset  {\mathcal S}_{\ell+1}$ of 15 multivariate Haar functions  with support 
$Q_{\ell,\alpha}$ given by all possible tensor products 
$$\eta_{\ell,\alpha,\beta}(x):=\tilde\eta_{\ell,\alpha_1,\beta_1}(x_1)\cdot
\tilde\eta_{\ell,\alpha_2,\beta_2}(x_2)\cdot\tilde\eta_{\ell,\alpha_3,\beta_3}(x_3)\cdot\tilde\eta_{\ell,\alpha_4,\beta_4}(x_4)$$
where 
$$\tilde\eta_{\ell,\alpha_k,\beta_k}(x_k):=
\begin{cases}
2^{\ell/2}\cdot 1_{[\alpha_k,\alpha_k+1)}(2^{\ell}x_k), \quad&\text{if }\beta_k=0\\
2^{\ell/2}\cdot\Big( 1_{[\alpha_k,\alpha_k+\frac12)}-
1_{[\alpha_k+\frac12,\alpha_k+1)}\Big)
(2^\ell x_k)
, \quad&\text{if }\beta_k=1
\end{cases}$$
for $k=1,2,3,4$.

For $\ell\in\N_0$, the block 
$${\mathcal H}_\ell:=
\{\eta_{\ell,\alpha,\beta}: \alpha\in A_\ell, \, \beta\in B_\ell\}
$$
consists of $15\cdot 2^{4\ell}$ Haar functions which we call Haar functions of level $\ell$. The union of all of them, 
$${\mathcal H}=\bigcup_{\ell=0}^\infty {\mathcal H}_\ell,$$ is a complete orthonormal system in $\mathring L^2(\TT^4)$, called \emph{isotropic 4-dimensional Haar system}, cf. \cite{Oswald}.
Moreover, 
\begin{equation}\label{discreteHaar}{\mathcal S}_{\ell}=\text{span}\bigg(\bigcup_{\kappa=0}^{\ell-1} {\mathcal H}_\kappa\bigg).
\end{equation}

For $x\in\TT^4$ and $\ell\in\N$, the unique cube $Q\in{\mathfrak Q}_\ell$ with $x\in Q$ will be denoted by $Q_\ell(x)$. Given a function $u\in \mathring L^1(\TT^4)$, we define the function $u_\ell\in{\mathcal S}_\ell$ by
\begin{equation}\label{proj-on-ql}
u_\ell(x):=2^{4\ell}\,\int_{Q_\ell(x)}u\,d\Leb^4.
\end{equation}
Restricted to $\mathring L^2(\TT^4)$, the map $\pi_{\mathfrak Q_\ell}: \, u\mapsto u_\ell$  is the $L^2$-projection onto the linear subspace ${\mathcal S}_\ell$. Moreover,
\begin{equation}\label{proj-on-eta}u_\ell=\sum_{\kappa=0}^{\ell-1} \sum_{\iota\in I_\kappa} \langle u,\eta_{\kappa,\iota}\rangle\, \eta_{\kappa,\iota}\,.
\end{equation}

\subsection{The Semi-discrete Gaussian Field}
Let an i.i.d. family of ${\mathcal N}(0,1)$ random variables $(\xi_{\ell,\iota})_{\ell\in\N_0, \iota\in I_\ell}$ be given with $I_\ell=A_\ell\times B_\ell$ as before.
For $\ell\in\N$ put
$$\hat h_\ell^\omega(x):= \sqrt8\,\pi\, \sum_{\kappa=0}^{\ell-1}\sum_{\iota\in I_\kappa} \xi^\omega_{\kappa,\iota}\cdot \mathring {\sf G}\,\eta_{\kappa,\iota}(x).$$
Here $\mathring{\sf G}$ denotes the grounded Green operator on the 4-torus, given as an integral operator
$\mathring{\sf G}u(x)=\int_{\TT^4} \mathring G(x,y) \,u(y)\,d\Leb^4(y)$
in terms of the grounded Green kernel $\mathring G(x,y)$ on $\TT^4$.
(For related results with the Green kernel of the torus replaced by the Green kernel of the discrete torus, see Section \ref{4.3} below.)
Moreover, define the non-symmetric kernel
$\mathring G_\ell(x,z):=\big(\pi_{\mathfrak Q_\ell}\mathring G(x,.)\big)(z)=2^{4\ell}\int_{Q_\ell(z)}\mathring G(x,v)\,d\Leb^4(v)$ and
put
$$\hat k_\ell(x,y):= 8\pi^2\,\int_{\TT^4}\mathring G_\ell(x,z)\mathring G_\ell(y,z)\,d\Leb^4(z).$$
As $\ell\to\infty$, this 
converges 
pointwise to 
$k(x,y):=8\pi^2\, \int_{\TT^4} \mathring G(x,z)\,\mathring G(z,y)\,d\Leb^4(z)$,
which --- up to the pre-factor --- is the Green kernel for the bi-Laplacian $\Delta^2$.

\begin{proposition} For every $\ell\in\N$, 
\begin{enumerate}

 \item for every $\omega$, the function $\hat h^\omega_\ell$
 \begin{itemize} 
\item is in $\C^1(\TT^4)$ and grounded (i.e. $\int_{\TT^4} \hat h_\ell^\omega\,d\Leb^4=0$);
\item is smooth off the boundaries of dyadic cubes $Q\in\mathfrak Q_\ell$;
\item has constant Laplacian on the interior of each dyadic cube $Q\in\mathfrak Q_\ell$;
\item is the sum $\sum_{\kappa=0}^{\ell-1}\sum_{\iota\in I_\kappa} \hat h^\omega_{\kappa,\iota}$ of functions $\hat h^\omega_{\kappa,\iota}= \sqrt8\,\pi\,  \xi^\omega_{\kappa,\iota}\cdot \mathring {\sf G}\,\eta_{\kappa,\iota}$ each of which is harmonic on  $\TT^4\setminus \bar Q_{\kappa,\iota}$ for the dyadic cube $Q_{\kappa,\iota}\in\mathfrak Q_\ell$;
\end{itemize}
\item for every $x\in \TT^4$, the random variable
 $\hat h_\ell(x)$ is centered 
 and  Gaussian  with variance $\hat k_\ell(x,x)$, the latter being independent of $x$;
 \item
 $\hat h_\ell$ is a centered Gaussian field with covariance function $\hat k_\ell(x,y)$.
\end{enumerate}
\end{proposition}

\begin{proof} We show c), the rest is straightforward. By the very definition of $\hat h_\ell$, the i.i.d.~property  of the $\xi_{\kappa,\iota}$ and the projection properties 
\eqref{proj-on-ql} \& \eqref{proj-on-eta}, 
\begin{align*}
{\mathbf E} \Big[\hat h_\ell(x)&\cdot \hat h_\ell(y)\Big]=8\pi^2\cdot\sum_{\kappa=0}^{\ell-1} \sum_{\iota\in I_\kappa}\big\langle \mathring{G}(x,.),\eta_{\kappa,\iota}\big\rangle\cdot \big\langle \mathring{G}(y,.),\eta_{\kappa,\iota}\big\rangle\\
&=8\pi^2\cdot\Big\langle \pi_{\mathfrak Q_\ell} \mathring{G}(x,.),\pi_{\mathfrak Q_\ell} \mathring{G}(y,.)\Big\rangle\\
&=8\pi^2\cdot 2^{8\ell}\,\int_{\TT^4}\bigg( \int_{Q_\ell(z)}\mathring G(x,v)\,d\Leb^4(v)\cdot\int_{Q_\ell(z)}\mathring G(y,w)\,d\Leb^4(w)\bigg)\,d\Leb^4(z)\\
&
=\hat k_\ell(x,y).
\end{align*}
\end{proof}

\begin{theorem} 
\begin{enumerate}

 \item The centered Gaussian random field $h$ with covariance function $k$ (as introduced and studied in section 3.1) is given in the case of the 4-torus   by
 $$h:=\sqrt 8 \, \pi\, \sum_{\kappa=0}^\infty\sum_{\iota\in I_\kappa} \xi_{\kappa,\iota}\cdot \mathring {\sf G}\eta_{\kappa,\iota}$$
 and called biharmonic Gaussian field.
\item The convergence $\hat h_\ell\to h$ holds in $L^2(\mathbf P)\times \Hil^{-\epsilon}(\TT^4)$ for every $\epsilon>0$. In particular,
for a.e.~$\omega$ and every $\epsilon>0$, 
$$h^\omega\in \Hil^{-\epsilon}(\TT^4).$$

\item For every $u\in\Hil^{-2}(\TT^4)$, the family $(\langle u,\hat h_\ell\rangle)_{\ell\in\N}$ is a centered $L^2(\mathbf P)$-bounded martingale and
$$\langle u,\hat h_\ell\rangle\to\langle u,h\rangle\quad\text{in $L^2(\mathbf P)$ as }\ell\to\infty.$$

 \end{enumerate}
\end{theorem}

\begin{proof}
a) For convergence (and well-definedness) of the infinite sum, see b) and/or c) below. To identify the covariance, observe that
\begin{align*}
{\mathbf E} \Big[\langle u,h\rangle^2\Big]&=
8\pi^2\cdot\sum_{\kappa=0}^\infty \sum_{\iota\in I_\kappa}\big\langle \mathring{G}u,\eta_{\kappa,\iota}\big\rangle^2
=8\pi^2\cdot\big\| \mathring{G}u\big\|^2=\big\langle u, ku\rangle.
\end{align*}

b) For $s>0$, let  $\mathring{\sf G}^{s}$ denote the $s$-th power of the operator $\mathring{\sf G}$ and let 
$\mathring G^{(s)}$ denote its kernel which is given by the formula
$$\mathring G^{(s)}(x,y)=\frac1{\Gamma(s)}\int_0^\infty t^{s-1}\mathring p_t(x,y)dt$$
in terms of the grounded heat kernel $\mathring p_t(x,y)=p_t(x,y)-1$.
Then for $\epsilon>0$,
\begin{align*}
\frac1{8\pi^2}\,{\mathbf E} \Big[\|h\|_{\Hil^{-\epsilon}}^2\Big]&=\frac1{8\pi^2}\,{\mathbf E} \Big[\big\| \mathring{\sf G}^{\epsilon} h\big\|_{L^2}^2\Big]
=\sum_{\kappa=0}^\infty \sum_{\iota\in I_\kappa} \big\| \mathring{\sf G}^{1+\epsilon} \eta_{\kappa,\iota}\big\|_{L^2}^2\\
&=\int_{\TT^4}\Big\|\mathring G^{(1+\epsilon)}(.,z)\Big\|^2\,d\Leb^4(z)=\int_{\TT^4}\mathring G^{(2+2\epsilon)}(z,z)\,d\Leb^4(z)\\
&=\mathring G^{(2+2\epsilon)}(0,0)<\infty
\end{align*}
since $\mathring G^{(s)}$ for $s>n/2=2$ is a bounded function, see \cite{DHKS}.
This proves that $\|h\|_{\Hil^{-\epsilon}}<\infty$ for a.e.~$\omega$. The convergence $\hat h_\ell\to h$ in $\Hil^{-\epsilon}(\TT^4)$ follows similarly.

c) By construction, for every $x\in\TT^4$, the family $(\hat h_\ell(x))_{\ell\in\N}$ is a centered martingale. The martingale property immediately carries over to the family $(\langle u,\hat h_\ell\rangle)_{\ell\in\N}$ for any function or distribution $u$.
The $L^2$-boundedness follows from
\begin{align*}
\frac1{8\pi^2}\cdot\sup_\ell {\mathbf E} \Big[\langle u,\hat h_\ell\rangle^2\Big]&=\sup_\ell \big\|\pi_{\mathfrak Q_\ell}\mathring{\sf G}u\big\|^2_{L^2}=\big\|\mathring{\sf G}u\big\|^2_{L^2}<\infty.
\end{align*}


\end{proof}


\begin{remark} a) With $\eta_{\ell,\iota}$ for $\iota=(\alpha,\beta)$ in $I_\ell$  as above and $\psi_{\ell,\iota}:=\mathring{\sf G}\eta_{\ell,\iota}$, the  family $(\psi_{\ell,\iota})_{\ell\in\N_0,\iota\in I_\ell}$ is a complete orthonormal system in the Hilbert space $\Hil^2(\TT^4)$, equipped with the scalar product
$$\langle u,v\rangle_{\Hil^2}:=\langle \Delta u,\Delta v\rangle_{L^2}.$$

b) Given any complete orthonormal system $(\psi_{k})_{k\in\N}$ in the Hilbert space $\mathring\Hil^2(\TT^4)$ and any  i.i.d. sequence of ${\mathcal N}(0,1)$ random variables $(\xi_{k})_{k\in\N}$, with the same arguments as for Theorem 8 one can prove that the Gaussian random field
$$h_\ell:= \sqrt8\,\pi\, \sum_{k=1}^{\ell} \xi_{k}\cdot \psi_{k}$$
converges in $L^2(\mathbf P)\times \Hil^{-\epsilon}(\TT^4)$ for every $\epsilon>0$ as $\ell\to\infty$ to the biharmonic Gaussian field $h$ on the 4-torus. Moreover, $\langle u,h_\ell\rangle \to \langle u,h\rangle$ in $L^2(\mathbf P)$ for every $u\in \Hil^{-2}(\TT^4)$.
\end{remark}
\subsection{The Semi-discrete Liouville Measure}\label{4.3}

Given a Gaussian field $h$ as considered in Theorem 8, then following Proposition \ref{const-approx-mart} a semi-discrete approximation of it is defined  by
$$h_\ell:=\pi_{\mathfrak Q_\ell} h\qquad (\forall \ell\in\N).$$
For each $\ell$, this is a centered Gaussian field with covariance given by
\begin{align*}
k_\ell(x,y)&:=8\pi^2\, 2^{8\ell}\int_{Q_\ell(x)}\int_{Q_\ell(y)}\bigg(\int_{\TT^4}\mathring G(z,v)\mathring G(z,w) \,d\Leb^4(z)
\bigg)d\Leb^4(w)\,d\Leb^4(v)\\
&=8\pi^2\, 2^{8\ell}\int_{Q_\ell(x)}\int_{Q_\ell(y)} k(v,w)\,d\Leb^4(w)\,d\Leb^4(v).
\end{align*}
%

For any $\gamma\in\R$ and $\ell\in\N$, we define the \emph{semi-discrete quantum Liouville measure} $\mu_\ell=\rho_\ell\,\Leb^4$ as the random measure on $\TT^4$ with density w.r.t.~$\Leb^4$ given by
$$\rho^\omega_\ell(x):=\ee^{\gamma h^\omega_\ell(x)-\frac{\gamma^2}2k_\ell(x,x)}.$$
\begin{corollary} For $|\gamma|<\sqrt 8$ and a.e.~$\omega$, the measures $\mu^\omega_\ell$ as $\ell\to\infty$ weakly converge to the Borel measure $\mu^\omega$   introduced and studied in section 3.2, $$\mu_\ell\to \mu\quad\mathbf P\text{-a.s.}$$
For $|\gamma|<2$, the convergence also holds in $L^2(\mathbf P)$.
\end{corollary}


\subsection{Discrete Random Objects}
To end up with fully discrete random objects, we have to replace the Green function on the continuous torus by the Green function on the discrete torus. To do so, for $\ell\in\N$ we define the \emph{discrete Laplacian} $\Delta_\ell$ acting on functions $u\in L^2(\TT^4_\ell)$ by
$$-\Delta_\ell u= 2^{2\ell+3}\cdot\big(u-p_\ell u\Big), \qquad p_\ell u(i):=\frac18\sum_{j\in J_\ell} u(i+j)$$ 
with $J_\ell:=\Big\{(k,0,0,0), (0,k,0,0), (0,0,k,0), (0,0,0,k): \ k\in\{-2^{-\ell}, 2^{-\ell}\}\Big\}$.
Note that the associated discrete Dirichlet form
$\mathcal E_\ell( u,u):=-\langle u, \Delta_\ell u\rangle_{L^2}$ on $L^2(\TT^4_\ell)$ has a positive spectral gap
$\lambda^1_\ell:=\inf\Big\{\frac{ \mathcal E_\ell( u,u)}{\|u\|^2_{L^2}}: \ u\in \mathring L^2(\TT^4_\ell)\Big\}.$
Furthermore,
we define the grounded transition kernel by $\dot p_\ell(i,j):=\frac18{\bf 1}_{J_\ell}(i-j)-2^{-4\ell}
$.

The \emph{discrete Green operator} acting on grounded functions $u\in \mathring L^2(\TT^4_\ell)$ is defined by
\begin{equation}
\dot{\sf G}_\ell u(i):=2^{-2\ell-3}\, \sum_{k=0}^\infty p_\ell^ku(i)=2^{-2\ell-3}\, \sum_{k=0}^\infty \dot p_\ell^ku(i)=\sum_{j\in \TT^4_\ell} \dot G_\ell(i,j) u(j)
\end{equation}
where $\dot G_\ell(i,j):=2^{-2\ell-3}\, \sum_{k=0}^\infty \dot p_\ell^k(i,j)$ denotes the grounded Green kernel.
The convergence of the operator sum is granted by the positivity of $\lambda_\ell^1$:
\begin{align*}
\big\| p_\ell\big\|_{\mathring L^2, \mathring L^2}\le
\big\| p_\ell^{1/2}\big\|_{\mathring L^2, \mathring L^2}^2&=\sup_{u\in \mathring L^2}\frac{ \langle u, p_\ell u\rangle}{\|u\|^2_{\mathring L^2}}=1-\inf_{u\in \mathring L^2}\frac{ \mathcal E_\ell( u,u)}{\|u\|^2_{\mathring L^2}}=1-\lambda_\ell^1<1.
\end{align*}
This operator is  then  extended to an operator acting on functions $u\in S_\ell$ by
$$\bar{\sf G}_\ell u(x):= \dot{\sf G}_\ell\Big( u|_{\TT_\ell^4}\Big)\big(2^{-\ell}\alpha\big)\qquad\quad\forall x\in Q_{\ell,\alpha}, \alpha\in A_\ell.$$
In terms of the (extended) discrete Green operator we define the \emph{discrete Gaussian field}
$$\dot h_\ell^\omega(i):= \sqrt8\,\pi\, \sum_{\kappa=0}^{\ell-1}\sum_{\iota\in I_\kappa} \xi^\omega_{\kappa,\iota}\cdot \dot {\sf G}_\ell\,\Big(\eta_{\kappa,\iota}|_{\TT_\ell^4}\Big)(i)$$
on $\TT^4_\ell$
and its piecewise constant extension 
$$\bar h_\ell^\omega(x):= \sqrt8\,\pi\, \sum_{\kappa=0}^{\ell-1}\sum_{\iota\in I_\kappa} \xi^\omega_{\kappa,\iota}\cdot \bar {\sf G}_\ell\,\eta_{\kappa,\iota}(x)$$
on $\TT^4$. For $\gamma\in\R$, the \emph{discrete quantum Liouville measure} is given by
$$\dot\mu_\ell^\omega:=2^{-4\ell}\, \exp\Big(-\frac{\gamma^2}2\dot k_\ell
\Big)\sum_{i\in \TT^4_\ell} \exp\Big(\gamma\dot h_\ell^\omega(i)
\Big)\, \delta_{i}
$$
with $$\dot k_\ell:=\mathbf E\big[\dot h_\ell(i)^2\big]=8\pi^2\cdot 2^{-4\ell}\,\sum_{j\in \TT^4_\ell}\dot{\sf G}_\ell(i,j)^2=\pi^2\cdot2^{-8\ell-3}\,\sum_{k=0}^\infty k\,\dot p_\ell^{k}(i,i),$$ independent of $i\in\TT^4_\ell$.

Alternatively --- and equivalent in distribution according to \eqref{discreteHaar} ---  we can define the discrete Gaussian field by
\begin{equation}\label{discField}\dot h_\ell^\omega(i):= \sqrt8\,\pi\,  \big(\dot {\sf G}_\ell\, \dot\xi^\omega\big) (i)= 2^{-2\ell-3/2}\,\pi\, \sum_{k=0}^\infty\sum_{j\in\TT_\ell^4} \dot p^k_\ell\,(i,j)\, \,\dot\xi^\omega_{j}\end{equation}
on $\TT^4_\ell$
with  a sequence of $\mathcal N(0,1)$-i.i.d.~random variables $(\xi_j)_{j\in\TT_\ell^4}$ and
$\dot\xi_i:=\xi_i-2^{-4\ell}\sum_{j\in\TT^4} \xi_{j}$. In other words,
$$-\Delta_\ell \dot h^\omega_\ell(i)= \sqrt8\,\pi\, \dot\xi^\omega_i.$$
Moreover,
$$\mathbf E\Big[\langle u,\dot h_\ell\rangle^2_{L^2}\Big]=8\pi^2\,\big\| \dot{\sf G}_\ell u\big\|^2_{L^2}\qquad \forall u\in L^2(\TT^4_\ell)$$
and thus the distribution of the Gaussian field $\dot h_\ell$  is explicitly given by the probability measure 
\begin{equation}
d\mathbf P_\ell(\zeta):=\frac1{Z_\ell}\, \exp\bigg(-\frac1{16\pi^2} \Big\| \Delta_\ell \zeta\Big\|^2_{L^2(\TT^4_\ell)}\bigg)\prod_{j\in\TT_\ell^4}\Leb^1(d\zeta_j),
\end{equation}
conditioned to the hyperplane $\big\{\sum_i\zeta_i=0\big\}$ 
in $\R^{\TT^4_\ell}$.
Here $\zeta=(\zeta_i)_{i\in\TT^4_\ell}$, $$\Big\| \Delta_\ell \zeta\Big\|^2_{L^2(\TT^4_\ell)}=2^{-4\ell}\sum_{i\in\TT^4_\ell} \Big|\zeta_i-\frac18\sum_{j\in J_\ell}\zeta_{i+j}\Big|^2
,$$
and $Z_\ell\in (0,\infty)$ denotes a suitable normalization constant.

The convergence $\bar h_\ell\to h$ and $\dot\mu_\ell\to\mu$ as $\ell\to\infty$ will be analyzed in detail in the forthcoming paper \cite{DHKS2}. For related convergence questions concerning biharmonic Gaussian random fields on the cube $[0,1]^4$ with Dirichlet boundary conditions, see \cite{Schweiger}.  
%
%
%
%

{\small
\bibliographystyle{plain}
\bibliography{confliouFuku.bib}

\begin{thebibliography}{10}

\bibitem{Cercle}
{Cercl{\'{e}}, B.}
\newblock {Liouville Conformal Field Theory on the higher-dimensional sphere}.
\newblock {\em arXiv:1912.09219}, 2019.

\bibitem{DHKS}
{Dello Schiavo, L.}, {Herry, R.}, {Kopfer, E.}, and {Sturm, K.-T.}
\newblock {Conformally invariant random fields, quantum Liouville measures, and
  random Paneitz operators on Riemannian manifolds of even dimension}.
\newblock {\em {Arxiv 2105.13925}}, 2021.

\bibitem{DHKS2}
{Dello Schiavo, L.}, {Herry, R.}, {Kopfer, E.}, and {Sturm, K.-T.}
\newblock {Polyharmonic fields and Liouville geometry in arbitrary dimension:
  from discrete to continuous}.
\newblock {\em {In preparation}}, 2021.

\bibitem{DDDF20}
{Ding, J.}, {Dub{\'{e}}dat, J.}, {Dunlap, A.}, and {Falconet, H.}
\newblock Tightness of {L}iouville first passage percolation for {$\gamma \in
  (0,2)$}.
\newblock {\em Publ.\ Math.\ Inst.\ Hautes \'{E}tudes Sci.}, 132:353--403,
  2020.

\bibitem{DRSV}
{Duplantier, B.}, {Rhodes, R.}, {Sheffield, S.}, and {Vargas, V.}
\newblock Log-correlated gaussian fields: an overview.
\newblock {\em Geometry, analysis and probability}, pages 191--216, 2017.

\bibitem{FOT}
{Fukushima, M.}, {Oshima, Y.}, and {Takeda, M.}
\newblock {\em {Dirichlet forms and symmetric Markov processes}}, volume~19 of
  {\em {De Gruyter Studies in Mathematics}}.
\newblock {de Gruyter}, extended edition, 2011.

\bibitem{GJMS}
{Graham, C.~R.}, {Jenne, R.}, {Mason, L.~J.}, and {Sparling, G.~A.~J.}
\newblock Conformally invariant powers of the {L}aplacian. {I}. {E}xistence.
\newblock {\em J.\ London Math.\ Soc.~(2)}, 46(3):557--565, 1992.

\bibitem{GM21}
{Gwynne, E.} and {Miller, J.}
\newblock Existence and uniqueness of the {L}iouville quantum gravity metric
  for {$\gamma\in(0,2)$}.
\newblock {\em Invent.\ Math.}, 223(1):213--333, 2021.

\bibitem{Kah85}
{Kahane, J.-P.}
\newblock {Sur le Chaos Multiplicatif}.
\newblock {\em {Ann.\ sc.\ math.\ Qu{\'{e}}bec}}, 9(2):105--150, 1985.

\bibitem{Lawler05}
{Lawler, G.~F.}
\newblock {\em Conformally invariant processes in the plane}, volume 114 of
  {\em Mathematical Surveys and Monographs}.
\newblock Amer.\ Math.\ Soc., Providence, RI, 2005.

\bibitem{Lawler18}
{Lawler, G.~F.}
\newblock Conformally invariant loop measures.
\newblock In {\em Proceedings of the {I}nternational {C}ongress of
  {M}athematicians---{R}io de {J}aneiro 2018. {V}ol. {I}. {P}lenary lectures},
  pages 669--703. World Sci.\ Publ., Hackensack, NJ, 2018.

\bibitem{LeGall19}
{Le Gall, J.-F.}
\newblock Brownian geometry.
\newblock {\em Jpn.\ J.~Math.}, 14(2):135--174, 2019.

\bibitem{LeGallMiermont12}
{Le Gall, J.-F.} and {Miermont, G.}
\newblock Scaling limits of random trees and planar maps.
\newblock In {\em Probability and statistical physics in two and more
  dimensions}, volume~15 of {\em Clay Math.\ Proc.}, pages 155--211. Amer.\
  Math.\ Soc., Providence, RI, 2012.

\bibitem{LSSW}
{Lodhia, A.}, {Sheffield, S.}, {Sun, X.}, and {Watson, S.S.}
\newblock {Fractional Gaussian fields: A survey}.
\newblock {\em {Probab.\ Surveys}}, 13:1--56, 2016.

\bibitem{MSImaginary1}
{Miller, J.} and {Sheffield, S.}
\newblock Imaginary geometry {I}: interacting {SLE}s.
\newblock {\em Probab.\ Theory Relat.\ Fields}, 164(3-4):553--705, 2016.

\bibitem{MSLQG1}
{Miller, J.} and {Sheffield, S.}
\newblock Liouville quantum gravity and the {B}rownian map {I}: the {${\rm
  QLE}(8/3,0)$} metric.
\newblock {\em Invent.\ Math.}, 219(1):75--152, 2020.

\bibitem{Oswald}
{Oswald:. P.}
\newblock {Haar system as Schauder basis in Besov spaces: The limiting cases
  for $0 < p \le 1$}.
\newblock {\em {INS-Preprint No. 1810, Bonn University}}.

\bibitem{Paneitz}
{Paneitz, S.~M.}
\newblock A quartic conformally covariant differential operator for arbitrary
  pseudo-{R}iemannian manifolds (summary).
\newblock {\em SIGMA Symmetry Integrability Geom. Methods Appl.}, 4:Paper 036,
  3, 1983.
\newblock Published in 2008.

\bibitem{Schramm07}
{Schramm, O.}
\newblock Conformally invariant scaling limits: an overview and a collection of
  problems.
\newblock In {\em International {C}ongress of {M}athematicians. {V}ol. {I}},
  pages 513--543. Eur. Math. Soc., Z\"{u}rich, 2007.

\bibitem{SchShe09}
{Schramm, O.} and {Sheffield, S.}
\newblock Contour lines of the two-dimensional discrete {G}aussian free field.
\newblock {\em Acta Math.}, 202(1):21--137, 2009.

\bibitem{SchShe13}
{Schramm, O.} and {Sheffield, S.}
\newblock A contour line of the continuum {G}aussian free field.
\newblock {\em Probab.\ Theory Relat.\ Fields}, 157(1-2):47--80, 2013.

\bibitem{Schweiger}
{Schweiger, F.}
\newblock {On the membrane model and the discrete Bilaplacian}.
\newblock {\em PhD Thesis Bonn University}, 2021.

\bibitem{She07}
{Sheffield, S.}
\newblock Gaussian free fields for mathematicians.
\newblock {\em Probab.\ Theory Relat.\ Fields}, 139(3-4):521--541, 2007.

\end{thebibliography}
}

\end{document}